\documentclass[11pt]{article}
\usepackage{amssymb, latexsym, mathtools, amsthm}
\begingroup
    \makeatletter
    \@for\theoremstyle:=definition,remark,plain\do{%
        \expandafter\g@addto@macro\csname th@\theoremstyle\endcsname{%
            \addtolength\thm@preskip\parskip
            }%
        }
\endgroup
\usepackage{enumerate}
\usepackage{pgf,tikz}
\usepackage{pdfsync}
\usepackage{txfonts}
\usepackage[T1]{fontenc}
\usepackage{booktabs}
 \usepackage{graphicx}
\definecolor{dnrbl}{rgb}{0,0,0.3}
\definecolor{dnrgr}{rgb}{0,0.3,0}
\definecolor{dnrre}{rgb}{0.5,0,0}
\usepackage[colorlinks=true, citecolor=dnrgr, linkcolor=dnrre, urlcolor=dnrre] {hyperref}
\usepackage{xcolor}
\usepackage{parskip}
\usepackage{tabularx, colortbl}
\usepackage{geometry}
\usepackage{etoolbox}
 \usepackage[scriptsize, up]{caption}
\usepackage{pgf,tikz}
\usetikzlibrary{decorations, decorations.pathmorphing}
\usetikzlibrary {shapes}
\theoremstyle{plain}

\theoremstyle{definition}

\makeatletter

\let\c@table\c@figure
\makeatother

\newcommand{\Nat}{\mathbb{N}}

\newcommand{\restr}{\upharpoonright}  
\newcommand{\un}{\uparrow} 
\newcommand{\de}{\downarrow} 

\newcommand{\bigo}[1]{\mathop{\bf O}\/\left({#1}\right)}

\newcommand{\DD}{\mathcal{D}}

\newcommand{\rcep}{$\mathbf{0}'$-right-c.e.\ }

\newcommand{\lcepp}{$\mathbf{0}^{(2)}$-left-c.e.\ }

\newcommand{\ml}{Martin-L\"{o}f }

\newcommand{\KG}{Ku\v{c}era-G{\'a}cs\ }
\newcommand{\pz}{$\Pi^0_1$\ }

\newcommand{\eg}{e.g.\ }
\newcommand{\ie}{i.e.\ }
\newcommand{\ce}{c.e.\ }

\newcommand{\lce}{left-c.e.\ }
\newcommand{\rce}{right-c.e.\ }

\newcommand{\pf}{prefix-free }

\renewenvironment{abstract}
 { \normalsize
  \list{}{
    \setlength{\leftmargin}{.0cm}%
    \setlength{\rightmargin}{\leftmargin}%
    }%
  \item {\bf \abstractname.} \relax}
 {\endlist}
\usepackage{titlesec}
 \titleformat{\paragraph}
{\normalfont\normalsize\em\bfseries}{\theparagraph}{1em}{$\blacktriangleright$\hspace{0.2cm}}
\titlespacing*{\paragraph}
{0pt}{3.25ex plus 1ex minus .2ex}{0.5ex plus .2ex}
\title{Aspects of Chaitin's Omega
\thanks{Barmpalias was supported by the 1000 Young Talents Program from the Chinese Government No. D1101130, NSFC grant No. 11750110425 and Grant No. ISCAS-2015-07 from the Institute of Software.
Thanks to Cris Calude, Rod Downey, Noam Greenberg, Denis Hirschfeldt, Kenshi Miyabe, 
Andr\'{e} Nies and Yu Liang, for valuable feedback on earlier versions of this article, and thanks
to Zhang Li and Fei Teng of the Chinese Academy of Sciences for their continual support.}}

\author{George Barmpalias}
\date{\today}
\begin{document}
\maketitle
\begin{abstract}
The halting probability of a Turing machine, also known as Chaitin's Omega, is an algorithmically random number with many interesting properties.
Since Chaitin's seminal work, many popular expositions have appeared, mainly focusing 
on the metamathematical or philosophical significance of Omega (or debating against it).
At the same time, a rich mathematical theory exploring the properties of Chaitin's Omega
has been brewing in various technical papers, which quietly reveals the significance
of this number to many aspects of contemporary algorithmic information theory.
The purpose of this survey is to expose these developments and tell a story
about Omega, which outlines its multifaceted mathematical properties and roles
in algorithmic randomness.
\end{abstract}
\vspace*{\fill}
\noindent{\bf George Barmpalias}\\[0.5em]
\noindent
State Key Lab of Computer Science, 
Institute of Software, Chinese Academy of Sciences, Beijing, China.\\[0.2em] 
\textit{E-mail:} \texttt{\textcolor{dnrgr}{barmpalias@gmail.com}}.
\textit{Web:} \texttt{\textcolor{dnrre}{http://barmpalias.net}}\par
\vfill \thispagestyle{empty}
\clearpage

\tableofcontents
\clearpage
\section{Introduction}
The two most influential contributions of Gregory Chaitin to the theory of algorithmic
information theory are (a) the information-theoretic extensions of G\"{o}del's incompleteness
theorem\footnote{Here
we mean G\"{o}del's first incompleteness theorem which asserts that every sufficiently 
powerful formal system is incomplete, in the sense that there are undecidable sentences
with respect to it. 
Chaitin's extensions of  
 G\"{o}del's incompleteness theorem were developed in 
 \cite{MR6822342,MR0429516, MR0455537,MR1183547}
and popularized in many articles and books including \cite{MR2381451,MR900163,MR1482519}.
We note that Kikuchi \cite{MR1477766} and Kritchman and Raz \cite{MR2761888}
have given proofs of the second incompleteness theorem in the spirit of Chaitin's arguments, using
Kolmogorov complexity.} and (b) the discovery of the halting probability
Omega (often denoted by $\Omega$), as a concrete algorithmically random real number.
Chaitin's numerous popular expositions of these discoveries have attracted
a certain amount of criticism, which concerns his philosophical interpretations
of the incompleteness results 
\cite{MR1026605,MR1663393} as well as his limited and subjective view of the
field of algorithmic information theory and its contributors \cite{gacsreview}.

The purpose of the present article is to expose a mathematical 
theory of halting probabilities which was developed in the last 40 years by numerous researchers
(much of it without Chaitin's active participation) and which reveals interesting and deep properties of
the number Omega. Most expositions on Omega in the literature
such as \cite{MR2371546,MR3444819,MR2089753,MR2604491} 
do very 
little\footnote{virtually nothing beyond 
the basic fact that it is algorithmically random and effectively 
approximable from below.} with regard to
its mathematical properties, either because they aim at a very general audience or because they
focus on its philosophical significance.
In contrast, we take the view that a better argument for the importance of
Chaitin's Omega -- one that is immune to attacks concerning its philosophical 
interpretations -- is a mathematical theory that reveals its complexities and its relevance
to many important topics in contemporary algorithmic information theory.

\subsection{About this survey}
We adopt a rather informal style of presentation, often omitting technical definitions
of established notions, either on the assumption that the reader is familiar with them or
on the basis that any ambiguity can be easily resolved by consulting a given technical reference.
At the cost of (lack of) self-containment, this approach will allow us to tell a coherent 
and concise mathematical 
story of Omega, based on a number of technical and relatively recent 
contributions that are absent from popular expositions.
The rich reference list and our multiple citations may be regarded as a compensation
for our approach to this survey. Due to the large number or theorems about Omega
that are discussed in the text, and in order to increase the readability of this article and
keep it concise, 
we have (a) avoided theorem displays and proofs, opting for a conversational mention of
the statements in the right context; and (b) 
suppressed the names of the
contributors of most results from the main text, merely citing the number of the
relevant bibliographic entry. 

The protagonist in this survey is Omega itself 
and the style of this presentation aims at
showing the impact of Chaitin's omega in contemporary research in algorithmic randomness
in the most straightforward way.
Our assumption about the reader is basic familiarity with Omega 
and the part of algorithmic information theory
which is sufficient to define it. 
For technical definitions of notions or exact statements of results that we
mention but did not include, the reader is
referred to the encyclopedic monographs by Downey and Hirschfeldt 
\cite{rodenisbook} and Li and Vitanyi \cite{MR1438307}, or the specific
research articles that are cited during the various discussions.
Inevitably, we do not aim or claim to be exhaustive; however we do strive to
include or at least mention most technical contributions about Omega that fit our narrative.
The reader will also find open problems, research suggestions and loose ends, 
in the context of many 
discussions about Omega.

%

\subsection{What is a halting probability?}
Suppose that we run a universal\footnote{By universal machine (in any standard model of Turing machines) we mean that it can effectively simulate any other Turing machine (in the class that it belongs to). In particular, given an effective number ing of all machines in the class $(M_e)$
there exists a computable function $e\mapsto \sigma_e$ from indices to strings such that 
$U(\sigma_e\ast\sigma)\simeq M_e(\sigma)$ for all $\sigma$. This  notion 
of universality also applies to \pf machines, and contrasts the notion of optimal \pf machines
in the context of Kolmogorov complexity (see \cite[Definition 2.1]{MR1438307}). Optimal
machines are the ones with respect to which the Kolmogorov complexity of any string
is minimal within an additive constant compared to any other machine.} Turing machine on a random (in the probabilistic sense) 
binary program. More precisely, whenever the next bit of the program is required during the
program execution, we flip a coin and feed the binary output to the machine.
On the basis of this thought experiment, which we are often going to refer to as
{\em Chaitin's thought experiment},  we can then consider the probability
that the universal Turing  machine $U$ will halt.
If we consider this problem in the context of self-delimiting machines\footnote{Self-delimiting machines, are Turing machines such that
their halting only depends on the length of the initial segment of the input that they
read, and not on the length of the given program. This restriction is only required for
the present expression of the halting probability in terms of finite programs. 
Alternatively one could consider a universal oracle Turing machine $M$
and define the halting probability as the uniform Lebesgue measure of the 
binary streams\footnote{By {\em binary stream} we mean an infinite binary sequence.} 
$X$ such that $M(X)$ halts on an empty program. 
Chaitin \cite{MR0411829} 
showed that self-delimiting machines are equivalent to Turing machines with 
\pf domain and vice-versa, in the sense that 
computations in one model can be simulated by computation in the other.
However Juedes and Lutz \cite{Juedes2000} showed that
but this simulation incurs an exponential blow-up on the running
time, unless $P=NP$.} then
the halting probability $\Omega_U$ of $U$ 
takes the following simple expression:
\[
\Omega_U=\sum_{U(\sigma)\de} 2^{-|\sigma|}
\]
where $\sigma$ represents the random finite binary programs
and $U(\sigma)\de$ denotes the fact that $U$ halts on input $\sigma$.
Note that what we call Omega or {\em halting probability} (or alternatively Omega number) 
is not a single real number,
but a family of real numbers indexed by the corresponding universal Turing 
machine. The relationships and differences of different Omega numbers 
will be explored in 
Section \ref{YfgHxGw4Ku}.
Also note that $\Omega_U$ is a \lce real, \ie it is the limit of an increasing computable
sequence of rationals.

One can view Omega as a compressed version of Turing's halting set
$H=\{\sigma\ |\ U(\sigma)\de\}$.
Indeed, if we coded the halting set into a binary stream $\eta$ in a canonical way,
ordering the strings first by length and then lexicographically, and representing
the outcome $U(\sigma)\de$ by 1 and the outcome $U(\sigma)\un$ by 0, then
then the first $n$ bits of $\Omega_U$ can give answers to the first $m$ bits of
$\eta$, where $m$ is exponentially larger than $n$. The precise relationship between
these numbers $n,m$ will be explored in Section \ref{yKBcWTkPDT}.

Any \pf machine $M$ has a halting probability $\Omega_M$ -- not just the universal ones.
All halting probabilities are \lce reals and vice-versa, as a consequence of the 
Kraft-Chaitin theorem, every \lce real is the halting probability of some \pf machine.
In \cite{presentcerials} it was shown that there are non-computable \lce reals $\alpha$
such that any \pf machine with halting probability $\alpha$ has a computable domain.

Returning to universal \pf machines, 
a fundamental result is 
\begin{equation}\label{nIfPvZ8kto}
\parbox{12cm}{{\bf The Omega characterization:} \ml random \lce reals are exactly the halting probabilities of universal \pf machines}
\end{equation}
and moreover the same holds if we replace `universal' with `optimal' in the sense of
Kolmogorov complexity. This was proved in
the cumulative work 
in \cite{Solovay:75,Calude.Hertling.ee:01,Kucera.Slaman:01}
and will be discussed in more depth in the main part of this survey.\footnote{There are a number of
self-contained presentations of this result including \cite{MR1872277,MR2965512}. In our view
the best presentation is in \cite[Chapter 9]{rodenisbook}.}
In \cite{tcs/CaludeNSS11} it was shown that there are non-universal and even non-optimal
\pf machines whose halting probability is 1-random, \ie the same as the halting probability of
a universal \pf machine.
Moreover it was shown that if the domain of a \pf machine $M$ contains the domain of a 
universal \pf machine, then the halting probability of $M$ is 1-random (though $M$ may not be
universal).

\subsection{The number of wisdom}\label{lNZ1nm3t28}
Just as the halting set holds the answers to many interesting mathematical problems
(Goldbach's Conjecture, Riemann's hypothesis or 
any one-quantifier definable problem in formal arithmetic) so does Omega. The difference is
that the first $n$ bits of Omega hold more answers than the first $n$ bits of
(the characteristic sequence of) $H$.
Some authors \cite{bennetgarden,MR2427555} 
have noted that the first 10,000 bits of Omega (with respect to a
canonical universal Turing machine) contain the answers to many open problems
in mathematics, although extracting such answers would take unrealistically long computations.
Calude et.al.\ \cite{measurecomrieman}
compute that the answer for Riemann's hypothesis is contained within the
first 7780 many bits of $\Omega_U$, for a canonical $U$ of their own design.\footnote{The same
canonical universal machine $U$ is used in the papers \cite{ijbc/CaludeD07,CaludeDS02}.}
Hence, in the words of Bennett \cite{bennettechrep}, this number 
{\em embodies an enormous amount of wisdom in a very small space}
and this property has been the main premise in popular articles such as Gardner \cite{bennetgarden}.

Chaitin \cite{MR1011475} showed how to
produce an exponential diophantine equation\footnote{an equation involving only addition, multiplication, and exponentiation of non-negative integer constants and variables} 
whose solution set defines Omega.\footnote{If we consider one variable as 
parameter, we obtain an infinite series of equations. Then
the $n$th bit of Omega is 0 or 1 according to whether
the $n$th equation has finitely or infinitely many non-negative integer solutions.} 
Variations of this expression of Omega were later obtained in  
\cite{ordkieu,Matiyasevich2005, MR2427549}.
In \cite{MR2484232} Chaitin also gives an expression of Omega in terms of
tilings of the half-plane and in \cite{MR2209802} he discusses an interpretation of the
approximation to Omega by cellular automata.
Mathematically, given the classic results which allow to represent computably enumerable
in terms of solutions  of diophantine equations or tilings of the plane, the above expressions
of Omega are not very surprising. In fact, as
G\'{a}cs \cite{gacsreview} points out, some of them are straightforward
consequences of known results and the fact that Omega is a \lce real. However for
Chaitin these examples serve a deeper purpose, namely 
a demonstration that a certain type of randomness exists in mathematics.
Such arguments about the importance of Omega 
can be found in his popular books such as \cite{MR3444819}.

\section{Omega in algorithmic information theory and metamathematics}
In this section we discuss some of the more basic properties of Chaitin's Omega.
In Section \ref{UR9pq9LKSa} we briefly expose its properties from the point of view of
algorithmic information theory. In Section \ref{feqbn5r8Fn} we discuss the other
popular aspect of Chaitin's Omega, which is its role in demonstrating undecidability
in formal systems.
In Section \ref{olMkNQSb75} we discuss how the halting probability 
relates to the concept of
algorithmic probability of Solomonoff, and how Chaitin's randomized machine thought experiment
can be used in order to express a number of algorithmic complexity properties
in terms of probabilities. This discussion is a step towards more advanced results
that will be discussed in the later sections of this survey.
In Section \ref{N5JNQAqApW} we show how algorithmic probability can be combined
with a recent result about differences of Omega numbers in order to derive
a surprising fact about the
probability of undecidable sentences in arithmetic, in the context of Chain's thought experiment.
In Section \ref{pqVFBlQdiE} we discuss Tadaki's `watered-down' versions of Omega
and in Section \ref{Rv7TFLwx6C}
we present an analogue of Omega for the class of computably enumerable sets.

\subsection{Randomness, Incompressibility and Unpredictability of Omega}\label{UR9pq9LKSa}
 The halting probability of a universal Turing machine
 is an interesting concept, but what is more interesting is its
 mathematical and algorithmic properties.
 Since the halting set is computably enumerable, it follows that
$\Omega_U$ is the limit of an increasing computable sequence of rationals -- it is a
{\em \lce real}. 
Chaitin  \cite{MR0411829} also showed that it is {\em algorithmically random} in the most
standard sense -- \ml random. 
Algorithmic randomness is a negative concept, characterized by {\em avoidance}
of statistical tests, {\em un}predictability with respect to effective betting strategies, or
{\em in}compressibility with respect to effective compression.

Intuitively, algorithmic randomness means lack of algorithmically identifiable
properties. In contrast, the fact that Omega is \lce is a positive property,
suggesting that this number is constructible in some algorithmic sense.
It is this contrast that makes Omega
special.\footnote{ Before Chaitin's discovery, the most concrete \ml random real
known was a 2-quantifier definable number exhibited in
Zvonkin and Levin \cite{MR0307889}.} 
Moreover, the halting probability has a specific mathematical meaning,
which contrasts the intuition that algorithmically random objects are 
unidentifiable. The Omega characterization \eqref{nIfPvZ8kto} showed that
these two opposing properties -- algorithmic randomness and computable enumerability
of the left Dedekind cut -- characterize the Omega numbers.

\subsection{Omega in the \lce reals}\label{LwCy1kABvY}
Solovay \cite{MR1793895} initiated the study of Omega as a member
of the class of \lce reals. In order to compare \lce reals according to how fast
they can be approximated by monotone rational sequences, he defined the
 {\em Solovay reducibility on \ce reals}, 
where $\alpha$ is Solovay reducible to $\beta$ if from any good 
rational approximation $q<\beta$ to $\beta$ we can effectively obtain a good approximation
$f(q)<\alpha$ to $\alpha$:
\[
\parbox{13cm}{{\bf Solovay reducibility} $\alpha<_S \beta$: there 
exists a partial computable function $f$ and a constant $c$ such that for each $q<\beta$
we have $f(q)\de <\alpha$ and $\alpha-f(q)< c\cdot (\beta-q)$.}
\]
Moreover he showed that Chaitin's Omega is of complete Solovay degree i.e., roughly
speaking, 
any good approximation to Omega encodes a good approximation to any given \lce real.
One of the by-products of the proof of the Omega characterization \eqref{nIfPvZ8kto}
was the converse of the latter statement, giving a 
characterization of Omega numbers as the complete \lce reals with respect to
Solovay reducibility, \ie the maximum elements with respect to this preorder. 

Coarser reducibilities measuring randomness were introduced in \cite{MR2030512},
giving further characterizations of Omega numbers as the complete \lce reals in
certain degree structures induced by the reducibilities.
For example, 
\begin{equation}\label{c5Ve5Gxpr}
\parbox{14cm}{a \lce real $\beta$ is an Omega number if and only if 
$K(\alpha\restr_n\ |\ \beta\restr_n)=\bigo{1}$ for all \lce $\alpha$}
\end{equation}
\ie if
each initial segment of any real $\alpha$ is basically coded in the corresponding initial segment
of $\beta$, modulo a fixed-length program. Finally,
\[
\parbox{14cm}{a \lce real $\beta$ is an Omega number if and only if its 
\pf initial segment complexity
dominates the initial segment \pf complexity of 
any other \lce real $\alpha$}
\]
modulo an additive constant, 
\ie if  $K(\alpha\restr_n)= K(\beta\restr_n)+\bigo{1}$ for all \lce reals $\alpha$.
Additional properties that characterize the Omega numbers inside the class of \lce reals
will be explored in Section \ref{8kFnCcyiiP}.

\subsection{Undecidability in formal systems and Omega}\label{feqbn5r8Fn}
Chaitin \cite{MR900163,MR1183547} 
expressed his incompleteness theorem in terms of Omega:
any formal system can determine only finitely many bits of Omega.\footnote{G\'{a}cs \cite{gacsreview} correctly points out that this is a consequence of Levin's classic work on
randomness and \lce semi-measures. A stronger version of
incompleteness in terms of Kolmogorov complexity is discussed 
in \cite[Section 2.7.1]{Li.Vitanyi:93}.}
He also combined this fact with the mathematical expressions of Omega (\eg in terms of
exponential diophantine equations) in order to argue again for the existence of randomness in
mathematics (see \cite{MR1482519} for many discussions on this). 
Solovay \cite{MR1793895} used a fixed-point construction in order to produce a
version of Omega for which  for which $ZFC$ cannot predict a
single bit. Calude \cite{MR1923902} combined Solovay's construction with
the Omega characterization in order to show that any Omega number,
is the halting probability of a certain universal \pf machine (provably in Peano Arithmetic)
such that $ZFC$ cannot prove any statement of the form `the $i$th bit of 
the binary expansion of Omega is $k$' for any position $i$ in the binary expansion of
Omega after the maximal prefix of 1s.

\subsection{Algorithmic probability and Omega}\label{olMkNQSb75}
Solomonoff \cite{Solomonoff:64} 
defined algorithmic probability 
in the context of inductive inference (see \cite[Chapter 4]{MR1438307}
for an up-to-date presentation).
The context here is that a universal machine gives a natural distribution of weight to
the various finite strings (which can be seen as codes for finite objects) which is often
called the a priori distribution. This is only a semi-measure (the total weight is less than 1)
and although it can be normalized into a probability distribution, as a semi-measure it
enjoys appealing effective properties: it is lower semi-computable or, in our terminology, \lce
as a function (it can be effectively approximated from the left).
Given any string $\sigma$, the a priori probability 
$P(\sigma)$ of $\sigma$ is the weight
of all the strings in the domain of the universal \pf machine which output $\sigma$, \ie
\[
P(\sigma)=\sum_{U(\rho)\de=\sigma} 2^{-|\rho|}
\]
Equivalently,  $P(\sigma)$ is the probability that the universal machine will output
$\sigma$, if supplied with random bits as a program.
Intuitively, the more mass $\sigma$ accumulates, the less complex it is expected to be. 
There is a simple but important result in algorithmic 
information theory which exemplifies this intuition:
 the \pf complexity $K(\sigma)$
of $\sigma$ equals (modulo an additive constant) 
the negative logarithm of its a priori probability 
$P(\sigma)$. Then the algorithmic probability $\Omega(A)$ of a set $A$ -- the probability that
the universal machine outputs a string in $A$ -- is simply the
sum of all $P(\sigma)$ for $\sigma\in A$ or equivalently:
\[
\Omega(A)=\sum_{U(\rho)\de\in A} 2^{-|\rho|}
\]
The properties of this number will be discussed in Section \ref{ITRpcDH23Q}.

In this context, many properties of strings or streams involving \pf complexity
can be expressed in terms of probabilities of events in the above thought experiment which
involves running the universal \pf machine on a random program.
For example, $K(n)$ (meaning the \pf complexity of the string consisting of $n$ many 0s) 
is, modulo an additive constant, the negative logarithm of the weight of all strings in the domain of $U$
such that $|U(\sigma)|=n$. Hence 
\begin{equation}
\parbox{14cm}{$K(n)$ is the negative logarithm of the probability that the randomized universal \pf machine outputs a string of length $n$.}
\end{equation}
One can also express conditional probability in terms of conditional \pf complexity. For example,
consider the property of a real $X$ that $K(X\restr_n\ |\ n^{\ast})=\bigo{1}$, where $n^{\ast}$ denotes
the shortest program for $n$, or equivalently the pair $(n, K(n))$.
Then for each $X$ we have
\begin{equation}\label{UWqzTh6inM}
\parbox{14cm}{$K(X\restr_n\ |\ n^{\ast})=\bigo{1}$ if and only if
the probability of obtaining output $X\restr_n$,
provided that the output has length $n$, has a positive lower bound.}
\end{equation}
Here is why: by symmetry of information we have 
$K(X\restr_n\ |\ n^{\ast})=K(X\restr_n)-K(n)$, while the probability 
of the property in \eqref{UWqzTh6inM}
is the probability of output $X\restr_n$ over the probability of an $n$-bit output, \ie
$2^{-K(X\restr_n)}/2^{-K(n)}$. Hence $K(X\restr_n\ |\ n^{\ast})=\bigo{1}$ means that
$K(X\restr_n)=K(n)+\bigo{1}$, which in turn is another way to say that 
$2^{-K(X\restr_n)}/2^{-K(n)}$ has a positive lower bound.
Incidentally, the property in \eqref{UWqzTh6inM} is no other than the well-known and
studied {\em $K$-triviality} (see \cite[Chapter 11]{rodenisbook} or \cite[Chapter 5]{Ottobook}).
Hence the  $K$-trivial streams are those whose initial segments are very likely to be produced
(in the above precise sense) when the universal \pf machine is run on a random input.

The study of $K$-triviality has been a significant part of 
research in algorithmic randomness in the last
15 years, and it is an area where Chaitin's Omega often plays an important role -- see
\cite{MR3248779} or the results discussed in Sections  
\ref{SSwIvk3kXz} and \ref{eGq9lCDoAE}.
 
\subsection{Algorithmic probability in formal systems and Omega}\label{N5JNQAqApW}
Let us now look at algorithmic probability in the context of formal systems.
The strings $\sigma$ which are described by the universal \pf machine may now be viewed
as sentences in the language of arithmetic.
Then the algorithmic probability of the provable sentences is some version of Omega, \ie
the halting probability of some other universal \pf machine. 
This is simply because the set of provable sentences is a \ce set $A$ and,
as observed by Chaitin \cite{chaitin2004algorithmic}, $\Omega(A)$ is \lce and 1-random
when $A$ is computably enumerable and nonempty.
A more surprising fact is the following:
\begin{equation}\label{I5BveuGRt3}
\parbox{14cm}{The algorithmic probability of the set of undecidable sentences in
formal arithmetic is \lce and 1-random.}
\end{equation}
This is curious: the undecidable sentences 
are not effectively verifiable or falsifiable, \ie
neither they nor their negations can be obtained by generating 
all sufficiently long proofs in arithmetic.
Yet their algorithmic probability behaves as
the algorithmic probability of the theorems of formal arithmetic, being 
effectively approximable from the left as if they where effectively enumerable.

This curious result is equivalent to saying that $\Omega_U(B)$ is \lce and 1-random
when $B$ is any non-empty effectively closed set. 
Note that $\Omega_U(B)$, just as the algorithmic probability of
the undecidable sentences in formal arithmetic, is a difference of \lce reals.
The proof of \eqref{I5BveuGRt3} 
appreared in \cite{omegax} and
relies heavily on the understanding of differences
of halting probabilities  and their approximations as these 
are discussed in Section \ref{xtMTpwLfAA}. On the top of these facts, the proof requires
what is known as a {\em decanter argument}, which is a sophisticated method 
mostly employed for the study of computationally weak sets, such as the $K$-trivial sets,
and was originally introduced in
\cite{MRtrivrealsH} for the establishment of the
incompleteness of  the $K$-trivial sets.

\subsection{A weakly random version of Omega}\label{pqVFBlQdiE}
Tadaki \cite{MR1888278} defined `watered-down' versions of Omega,
in the sense that they are less random and have faster approximations, as
\[
\Omega^s=\sum_{U(\sigma)\de} 2^{-\frac{|\sigma|}{s}}
\hspace{0.5cm}\textrm{for $s\in (0,1]$}
\]
and he showed that $\Omega^s$ is {\em weakly $s$-random}
in the sense that
\begin{equation}\label{KiLbBuS9iB}
\exists c\ \forall n\ K(\Omega^s\restr_n)\geq s\cdot n-c
\end{equation}
Equivalent conditions
of \eqref{KiLbBuS9iB} in terms of \ml tests or effective betting strategies,
just as in the usual notion of \ml randomness, can be found in 
\cite[Section 13.5]{rodenisbook}, where 
a direct relation to effective Hausdorff dimension is also established.

Note that when $s<1$, the number $\Omega^s$ is less compressed than Omega because
each convergent computation adds a smaller amount to the probability.
In \cite{MR2539532} the analogue of the Omega characterization 
\eqref{nIfPvZ8kto} 
was carried out for these weaker versions of Omega, including the characterization in terms
of speed of convergence (mainly the results 
in \cite{Solovay:75,Calude.Hertling.ee:01,Kucera.Slaman:01} which are
discussed in Section \ref{8kFnCcyiiP} of the present article). Further work on weak versions
of Omega can be found in \cite{MR2799279} while a different version of halting probability
is introduced in \cite{MR2269901}.

\subsection{Analogues of Omega in the computably enumerable sets}\label{Rv7TFLwx6C}
Is there an analogue of Omega in the \ce sets? Here we are  not necessarily looking for
a probability, but a real which is computably enumerable as a set, and whose initial segments
are universal or maximally complex in some sense akin to the properties of Omega explored
in Section \ref{LwCy1kABvY}.
We can either look at the initial segment complexity of the \ce sets and ask that it is sufficiently high,
or look at a reducibility amongst \ce sets that measures complexity, and consider the complete
sets with respect to this reducibility (if they exist). Our answer will satisfy both of these heuristics,
and this fact will support our bid for an analogue of Omega in the \ce sets.

Considering the reducibilities of Section \ref{LwCy1kABvY}, we start with
Solovay reducibility on the \ce sets, which measures 
{\em hardness of approximation}.
Hence we could look for the \ce set analogue of Omega in the class of \ce sets
which are the hardest to approximate in this context.
Unfortunately, it was discovered in \cite{cie/Barmpalias05} that 
there is no complete \ce set in the Solovay degrees  
and, even worse, for each \ce set $A$ there exists a \ce set $B$ of strictly larger
Solovay degree than $A$.
Hence we need to consider a coarser reducibility. 

Consider the reducibility implicit in \eqref{c5Ve5Gxpr},
which was introduced in \cite{MR2030512} by the name of
{relative $K$ reducibility}: $X\leq_{rK} Y$ when $K(X\restr_n\ |\ Y\restr_n)=\bigo{1}$;
moreover this is equivalent to $C(X\restr_n\ |\ Y\restr_n)=\bigo{1}$.
We have seen that amongst the \lce reals, Omega numbers are characterized
as the $\leq_{rK}$-complete \lce reals. 
Is there
an $rK$-complete \ce set, \ie a \ce set $A$ such that $K(W\restr_n\ |\ A\restr_n)=\bigo{1}$
for all \ce sets $W$? Surprisingly, the answer is yes by \cite{ipl/BarmpaliasHLM13}, 
and moreover the following are equivalent
for any \ce set $A$:
\begin{enumerate}[\hspace{0.5cm}(a)]
\item $A$ is Turing complete with respect to a linear oracle-use function;\footnote{One can define $X\leq_{lin} Y$ if $X$ is
Turing computable from $Y$ with linear oracle-use $n\mapsto an+b$ (or equivalently
$n\mapsto an$) for some positive constants $a,b$. Then $A$ is linearly complete in the \ce
sets if $W\leq_{lin} A$ for all \ce sets $W$.}
\item  $K(W\restr_n\ |\ A\restr_n)=\bigo{1}$ or
equivalently  $C(W\restr_n\ |\ A\restr_n)=\bigo{1}$ for all \ce sets $W$;
\item there exists $c$ such that  $C(A\restr_n)\geq \log n-c$ for all $n$.
\end{enumerate}
Given that every \ce set $W$ has a constant $c$ and  
infinitely many $n$ with $C(W\restr_n)\leq \log n+c$ this is arguably the best possible candidate
as an analogue of Omega in the \ce sets in terms of initial segment complexity.\footnote{Here 
we also note that the plain initial
segment complexity of any \ce set is bounded above by $2\log n+\bigo{1}$ and this is optimal up
to an additive constant \cite{BarzdinsCe,siamcomp/Kummer96}. With respect to \pf complexity
the upper bound is $(2+\epsilon)\cdot\log n$ for any $\epsilon>0$ and it sometimes fails
for $\epsilon=0$ by \cite{koba_rod}.}
Moreover there are natural sets in this class, namely the halting sets of universal machines with
respect to Kolmogorov numberings.\footnote{Kolmogorov numbering is a G\"{o}del numbering
to which every other computable
numbering can be reduced via a linearly bounded function.}
Hence 
\[
\parbox{15cm}{the halting sets with respect to Kolmogorov numberings 
are analogues of Omega in the \ce sets.}
\]
There are further connections between such canonical halting sets and Chaitin's Omega. 
Consider a computable order  $g$  (\ie nondecreasing and unbounded) such
that $\sum_i 2^{-g(i)}$ is
an Omega number. 
Such functions were studied in \cite{DBLP:conf/stacs/BienvenuD09,BienvenuMN11,Bienvenu11MNjou}
by the name of {\em Solovay functions} and, given that
Omega numbers have the slowest approximations amongst the \lce reals, 
 they are
{\em very slow-growing}.
In \cite{MR3521997} it was shown that 
\[
\parbox{14cm}{a \lce real is an Omega number if and only
if it computes a halting set with respect to a Kolmogorov numbering
with oracle-use a non-decreasing Solovay function.}
\]
In the same paper it was shown that every 
halting set with respect to a Kolmogorov numbering
computes Omega  with use $\bigo{2^{n}}$.
The reader may compare these results with Tadaki \cite{Tadaki:200972}
who considered halting sets $W$ of universal \pf machines, instead of the
more compact halting sets of universal plain Turing machines. It was shown that in the 
\pf case, about
$2^{n+2\log n}$ bits of $W$ are needed for the computation of $\Omega\restr_n$,
and that $2^{n+\log n}$ do not always suffice. On the other direction, the first $n$ bits of
Omega compute the first $2^n$ bits of $W$, and this is optimal up to a multiplicative constant
on $2^n$.

\section{Computable approximations to Omega}\label{8kFnCcyiiP}
The Omega characterization \eqref{nIfPvZ8kto} relied on the study of the 
monotone computable approximations
properties of Omega, which was initiated by 
Solovay in \cite{Solovay:75}
and was briefly discussed
in Section \ref{LwCy1kABvY}.
In the present section we discuss the approximation to Omega more thoroughly,
including some recent results.

\subsection{Rates of convergence amongst \lce reals}\label{ZGwomJvEL8}
Recall the Solovay reducibility from Section \ref{LwCy1kABvY}
which was used in order to measure
the speed of \lce approximations to \lce reals.  
Solovay showed that the
induced degree
structure, known as the Solovay degrees, 
has a maximum element and he called the members of the maximum degree 
$\Omega$-like reals. He also showed that any universal halting probability is $\Omega$-like.
Then the work in \cite{Calude.Hertling.ee:01,Kucera.Slaman:01} showed that
the $\Omega$-like reals are exactly the halting probabilities of universal \pf 
machines.\footnote{The same argument shows that this equivalence is not sensitive to
many features of the universal machine. For example, the equivalence holds if
we consider a universal oracle Turing machine and define its halting probability
as the measure of oracles $X$ which make the machine halt with input the empty string.
The same can be said of other models such as the monotone machines of 
Levin \cite{levinthesis,Levin:73} and even the optimal machines \cite[Definition 2.1]{MR1438307} in
the sense of Kolmogorov complexity. However when considering
the halting probability restricted to a \pz set of outputs, the randomness of
this real is robust for universal machines by \cite{omegax} but may not be robust for
mere optimal machines by \cite{jc/FigueiraSW06}.}

The proof that the omega numbers are exactly the halting probabilities of universal \pf machines
is, in a sense, non-uniform.
Given a \lce index of a 1-random \lce real $\alpha$ in $(0,1)$ the known argument obtains
an index of a universal \pf machine $U$ such that $\Omega_U=\alpha$ 
by non-effective means.\footnote{One constructs a certain \ml test $(V_i)$ and uses
a number $n$ such that  $\alpha\not\in V_i$ for all $i>n$ in order to define
a universal \pf machine $U$ such that $\Omega_U=\alpha$. The main source of non-uniformity is
the choice of this number $n$.}
The question whether this non-uniformity is necessary, is open.
In fact, many results that touch on this characterization of omega numbers are proved by 
non-uniform arguments, though the necessity of this non-uniformity has not been
established. For example, in \cite{lata} it was shown that for each universal \pf machine $U$
there exists another universal \pf machine $V$ such that $\Omega_U\neq \Omega_V+\beta$
for every \lce real -- the proof was non-uniform and the necessity of the non-uniformity was
left open.

The Solovay reducibility $\beta\leq_S\alpha$ 
between \lce reals $\alpha,\beta$ can be defined equivalently by
any of the following clauses:
\begin{enumerate}[\hspace{1cm}(a)]
\item there exists a rational $q$ such that $q\alpha-\beta$ is \lce
\item there exist a rational $q$ and $(\alpha_s)\to\alpha$, $(\beta_s)\to\beta$ such that 
$\beta-\beta_s< q\cdot (\alpha-\alpha_s)$ for all $s$;
\item there exist a rational $q$ and $(\alpha_s)\to\alpha$, $(\beta_s)\to\beta$ such that 
$\beta_{s+1}-\beta_s< q\cdot (\alpha_{s+1}-\alpha_s)$ for all $s$.
\end{enumerate}
Note that the set of rationals $q$ for which one of the above clauses holds is upward closed - if the clause holds for the rational  $q$ then it also holds for all rationals $q'>q$.  
Although it is not explicitly stated in \cite{Downey02randomness}, it follows from the proofs
that when $\beta \leq_S \alpha$,  the infimums of the rationals $q$ for which the clauses (a), (b) and (c) hold are equal.
A thorough study of the algebraic
aspects of the structure of Solovay degrees of \lce reals was undertaken in \cite{Downey02randomness,MR2320380}.

There are many ways in which any approximation to Omega is very much slower
than any approximation to any \lce real which is not an Omega number.
For example, it was shown  in \cite{Downey02randomness} (also see \cite{lata}) that 
if $(\alpha_s)$, $(\Omega_s)$  are any \lce approximations
to $\alpha,\Omega$ respectively  then
\begin{equation}\label{k2SthHLQp2}
\textrm{$\alpha$ is not an Omega number}\ \ \ \iff\ \ \ 
\lim_s \frac{\alpha-\alpha_s}{\Omega-\Omega_s}=0.
\end{equation}
This shows that Omega is much more intractable than any \lce real which is not an
Omega number. 

\subsection{Initial segment complexity of Omega numbers}
The intuition that Omega is considerably more intractable than any \lce real
which is not an Omega number is also reflected in a result of \cite{Demuth:75ce} 
(also see \cite{Downey02randomness} and 
\cite[Remark 3.5]{Kucera.Slaman:01}) which says that 
\begin{equation}\label{vYpY1QoUT}
\parbox{14cm}{Omega cannot be written as the sum
of two \lce reals which are not Omega numbers.}
\end{equation}
We note that by \cite{MRtrivrealsH}, 
for any \lce reals $\alpha,\beta$, the \pf complexity of
$\alpha+\beta$ is (within a constant) the maximum
of the \pf complexities of $\alpha$ and $\beta$.
Hence, using the definition of randomness in terms of 
\pf complexity, another way to write \eqref{vYpY1QoUT} is:
\begin{equation*}
\parbox{14cm}{if $\alpha,\beta$ are \lce and not random then 
$\limsup_n \Big(n-\max\{K(\alpha\restr_n), K(\beta\restr_n)\}\Big)=\infty$.}
\end{equation*}
In other words, if
$\alpha,\beta$ are \lce and  not Omega numbers, their complexity at position $n$ 
simultaneously drops 
well below $n$ for infinitely many $n$. It is interesting that by \cite{MR3248779}, 
a similar property
occurs with regard to the \lce reals which are not $K$-trivial:
\begin{equation*}
\parbox{14cm}{if $\alpha,\beta$ are \lce and not $K$-trivial then 
$\limsup_n \Big(\min\{K(\alpha\restr_n), K(\beta\restr_n)\}-K(n)\Big)=\infty$.}
\end{equation*}
In other words, if two \lce reals have non-trivial \pf complexity, then there
are infinitely many lengths $n$ at which their complexity simultaneously rises well above
$K(n)$.

It is tempting to seek a strengthening of \eqref{vYpY1QoUT} 
in the statement that initial segment \pf complexity of any non-random \lce
real diverges from the initial segment \pf complexity of Omega.\footnote{This question
was asked by Kenshi Miyabe during
the conference {\em Aspects of Computation} organized by the Institute of Mathematical
Sciences in Singapore in August and September 2017, where Yu Liang pointed to
the negative answer we describe below.}
Such a result would reinforce the intuition that Omega is much more complex than any 
non-random \lce real. Unfortunately,
results from \cite{MR2030512,jflBaiasL06} show that this is not true:
\begin{equation}\label{OUrdYDo5xf}
\parbox{14cm}{There exists a \lce real $\alpha$ which is not 1-random but
$\liminf_n \Big(K(\Omega\restr_n)-K(\alpha\restr_{n})\Big)<\infty$.}
\end{equation}
Indeed, in \cite{jflBaiasL06} it was shown that there are \lce reals which are
not computed by any Omega number with oracle-use $n\mapsto n+\bigo{1}$.
Then \eqref{OUrdYDo5xf} follows from the above result combined with the theorem in
\cite{MR2030512} that
if  $\lim_n (K(\Omega\restr_n)-K(\alpha\restr_{n}))$ is infinite and $\alpha$ is a \lce real
then $\alpha$ is computable from $\Omega$ with oracle-use 
$n\mapsto n+\bigo{1}$. Moreover, in \cite{IOPORT.05678491} it was shown that
a \ce degree contains a real as the one constructed in  \cite{jflBaiasL06} if and only if
it is array non-computable. Hence by the same argument we get the `only if' 
direction of the following:
\begin{equation*}
\parbox{14cm}{A \ce degree is array non-computable if and only if it contains a \lce
real $\alpha$ such that
$\liminf_n \Big(K(\Omega\restr_n)-K(\alpha\restr_{n})\Big)<\infty$.}
\end{equation*}
The `if' direction of this statement follows from the fact that the \pf initial segment
complexity of a \lce real of array computable degree is well below $n$.\footnote{In fact, this complexity
is bounded above by $2\log n$, \eg by \cite[Theorem 8.2.29]{Ottobook}, the fact that facile reals
have \pf complexity bounded above by $K(n)+g(n)+\bigo{1}$ for any computable order $g$, and the
fact that in the \ce degrees the array computable degrees coincide with the \ce traceable degrees.}

Solovay \cite{Solovay:75} produced a number of technical results regarding
the initial segment complexity of Omega, which are reproduced in
\cite[Section 10.2]{rodenisbook}.
A consequence of these results is that the initial segment \pf complexity
of Omega does not bound above the initial segment \pf complexity of any 2-random real, 
\ie any real which is random relative to the halting problem.
In fact, we can elaborate on this statement. Since $\Omega$ is 1-random, we have that
$K(\Omega\restr_n)$ is in the interval $(n,n+K(n))$ (modulo an additive constant). By Chaitin \cite{chaitinincomp} we have
\[
\lim_n \Big(K(\Omega\restr_n)-n\Big)=\infty
\]
and this, in fact, is a characterizing property of every 1-random real number.
On the other hand, Miller \cite{joelowklowwea} showed that
$\liminf_n \big(n+K(n)-K(X\restr_n)\big)<\infty$
is equivalent to the property that $X$ is 2-random, hence
\[
\liminf_n \Big(n+K(n)-K(\Omega\restr_n)\Big)=\infty.
\]
These results highlight the fact that Omega is not as random as a typical algorithmically
random number, a theme that is explored further in Section \ref{eGq9lCDoAE}.

\subsection{Comparing the rates of convergence amongst Omega numbers}\label{xtMTpwLfAA}
There is one aspect of approximations to Omega that remained unresolved until very recently,
where it was settled in \cite{omegax}. We have seen in \eqref{LHkfIXzYnx} that Omega numbers
are in their own league as far as the rate of their monotone effective approximation is concerned. 
But how do the approximations to two different Omega numbers related?
We know by \eqref{k2SthHLQp2} 
that they are both very slow compared to any approximation to any real
which is not an Omega number, but can we compare the two rates with each other?
From the Omega characterization, and in particular from \cite{Kucera.Slaman:01}, 
we know that given any \lce approximations
$(\omega_s)\to\omega$, $(\Omega_s\to \Omega)$ to two Omega numbers, we have
$\liminf_s \big[(\Omega-\Omega_s)/(\omega-\omega_s)\big]>0$.\footnote{This holds even if
$\omega$ is not an Omega number.} But does this limit exist?

Quite remarkably, in \cite{omegax} it was shown that 
given any \lce approximations
$(\alpha_s)\to\alpha$, $(\beta_s)\to \beta$,
\begin{equation}\label{LHkfIXzYnx}
\textrm{\small $\alpha,\beta$ are Omega numbers}\ \ \ \Rightarrow\ \ \ 
\DD(\alpha,\beta):=\lim_s \frac{\alpha-\alpha_s}{\beta-\beta_s}\ \ \ \ \ \parbox{6.5cm}{\small  exists and is positive and independent of the chosen approximations $(\alpha_s)\to\alpha$, $(\beta_s)\to \beta$.}
\end{equation}
Moreover this limit $\DD(\alpha,\beta)$ has a rather special property -- it is the
\begin{itemize}
\item infimum of all rationals $q$
with the property that $q\cdot\beta-\alpha$ is a \lce real;
\item supremum of all rationals $p$ such that 
$p\cdot\beta-\alpha$ is a \rce real.
\end{itemize}

In fact, given two Omega numbers $\alpha<\beta$, the value of $\DD(\alpha,\beta)$
determines whether $\beta-\alpha$ is an Omega number:
\begin{itemize}
\item if  $\DD(\alpha,\beta)<1$ then $\beta-\alpha$ is an Omega number;
\item if  $\DD(\alpha,\beta)\geq 1$ then $\beta-\alpha$ is not an Omega number;
\end{itemize}
More specifically, we have the following trichotomy:
\begin{itemize}
\item if  $\DD(\alpha,\beta)>1$ then $\beta-\alpha$ is \rce  and 1-random;
\item if  $\DD(\alpha,\beta)<1$ then $\beta-\alpha$ is \lce  and 1-random;
\item if  $\DD(\alpha,\beta)=1$ then $\beta-\alpha$ is not \rce or \lce or 1-random;
\end{itemize}

Based on \eqref{LHkfIXzYnx}, a
theory of derivation for the field of differences of \lce reals
was developed in \cite{derivationmiller}.

\subsection{Speeding-up the approximation to Omega}\label{QgvOh5fZzw}
We are interested in functions $f$ which grow sufficiently fast so that
we can use them to approximate Omega faster than the universal effective approximation.
This means that 
\begin{equation}\label{hQRfkTbkrD}
\textrm{for each $c$ there exists $n$ such that \hspace{0.2cm}
$\Omega-\Omega_{f(n)}< 2^{-c}\cdot (\Omega-\Omega_{n})$.}
\end{equation}
It is not hard to see that no computable function $f$ has this property.
In fact, it was shown in \cite{MR3248779} that if $f$ is low for $\Omega$
(\ie it is computable from a set relative to which $\Omega$ is random) then it does not have
the above property. It was also shown that every \ce degree which is not $K$-trivial
contains a function $f$ such that \eqref{hQRfkTbkrD} holds. Since the low for $\Omega$
\ce degrees are exactly the $K$-trivial \ce degrees, it follows that in the \ce degrees
exactly the ones that are not $K$-trivial compute (or even contain) functions $f$ with the
property \eqref{hQRfkTbkrD}.
\begin{equation}\label{FnBRrTNJWJ}
\textrm{A \ce degree can speed-up the approximation to Omega if and only if it is not $K$-trivial.}
\end{equation}
Moreover, if a \ce set $A$ is not $K$-trivial then there exists a canonical function in the degree
of $A$ with the
property \eqref{hQRfkTbkrD}, namely the settling time in any computable enumeration of $A$.

This is interesting! Combined with known characterizations of the low for $\Omega$
sets (see Section \ref{SSwIvk3kXz}), we have that
the following are equivalent for \ce sets:
\begin{enumerate}[\hspace{0.5cm}(a)]
\item $A$ can speed-up the approximation to Omega;
\item $\Omega^A_U$ is not a \lce real;
\end{enumerate}
where $U$ is any universal \pf machine.

Can the characterization \eqref{FnBRrTNJWJ} 
be generalized outside the class of \ce degrees?
A reasonable guess, given the facts from \cite{MR3248779},  
would be to conjecture that a degree computes $f$
with the
property \eqref{hQRfkTbkrD} if and only if it is not low for $\Omega$.
Although one direction of this equivalence is true, surprisingly the other is not.
Miller and Nies, see \cite[Section 8.1]{Ottobook}, showed that
non-computable low for $\Omega$ sets are necessarily hyperimmune, \ie 
they compute functions which are not
dominated by any computable function. This means that non-computable hyperimmune-free degrees
are not low for $\Omega$, and these clearly do not compute functions with the property  \eqref{hQRfkTbkrD}, since every function they compute a is dominated by a computable function.

On the other hand, given that there exists a function $f\leq_{tt} \emptyset'$
such that $\Omega-\Omega_{f(n)}< 2^{-n}\cdot (\Omega-\Omega_n)$, 
it is not hard to see that every array non-computable degree can 
speed-up the approximation to Omega.
The array computable degrees contain the low for Omega degrees, but also contain other degrees
which can speed-up the approximation to Omega (such as \ce degrees which are not $K$-trivial).
The question of exactly which degrees can speed up the approximation to Omega remains open.

\section{Omega and computable enumerability}\label{yKBcWTkPDT}
We have discussed that Omega can be seen as a compressed version of the halting problem,
and moreover it enjoys several completeness properties with respect to the \lce reals.
In this section we present a number of results that make the connection
between Omega and \ce sets or other \lce reals precise.
In Section \ref{8VvCMdWp4N} we discuss the problem of how many bits of Omega
are needed for the computation of $n$ bits of a \ce set like the halting problem, or a \lce 
real, and vice-versa
how many bits of a \ce set are needed in order to compute $n$ bits of Omega.
In Section \ref{YfgHxGw4Ku} we give a thorough examination to the question of
how similar and how different two universal halting probabilities can be.

\subsection{Omega and halting problems}\label{8VvCMdWp4N}
We have seen that many open problems in mathematics have their solutions coded in the
first few thousands bits of a canonical version of Omega, while some authors have worked-out more
precise bounds of this type. Technically speaking, given that such solutions are merely
answers to certain halting problems, the issue here is which initial segment of Omega is sufficient
to answer a given halting problem. It is possible to obtain asymptotics that show the 
lengths of the initial segments of Omega in relation to the amount of halting problems that 
they encode.
The first results of this type were obtained by Solovay \cite{Solovay:75} and first published in
\cite[Section 3.13]{rodenisbook}.
If we define
\begin{eqnarray*}
p(n)&=|\{\sigma\in 2^n\ |\ U(\sigma)\de\}|&\\
\DD_n&=\{\sigma\in 2^{\leq n}\ |\ U(\sigma)\de\} &\textrm{and}\ \ \  P(n)=|\DD_n|
\end{eqnarray*}
then 
$p(n)\sim P(n) \sim 2^{n-K(n)}$ and $\DD_n$ is uniformly computable from $\Omega\restr_n$.
Moreover $K(\Omega\restr_n\ |\ \DD_{n+K(n)})=\bigo{1}$, which means that the first 
$n+K(n)$ bits of the halting problem can be used for the computation of the first $n$ bits of 
Omega (with an additional program of fixed length which might not be uniformly given in $n$).

Tadaki \cite{Tadaki:200972} improved on the latter result by showing the following 
characterization.\footnote{Formally, 
there exists an oracle Turing machine $M$ such that
for each $n$ we have $M(\DD_{n+f(n)+\bigo{1}}(W))(n)=\Omega_V\restr_n$.}
\begin{equation}\label{1DeccmxzPv}
\parbox{14cm}{Given 
optimal machines $W$, $V$ and a computable function $f$, 
we have $\sum_i 2^{-f(i)}<\infty$ if and only if
$\DD_{n}(W)$ uniformly computes the first $n-f(n)-c$ bits of 
$\Omega_V$ for some constant $c$.}
\end{equation}
Here we may use
$2\log n$ as a simple representative of the functions 
$f$ with the property
$\sum_i 2^{-f(i)}<\infty$.
Then \eqref{1DeccmxzPv} says that 
given any optimal machines $W,V$, in order to compute $n-2\log n$ bits 
of the halting probability with respect to $V$, we need to know the halting problem 
with respect to $W$, for all strings of length at most $n$.
 
Tadaki also gave an analogous result concerning computations of halting sets from
halting probabilities.
\begin{equation*}
\parbox{14cm}{Given optimal machines 
$W$, $V$ and a computable function $f$, 
we have $f=\bigo{1}$ if and only if the first $n$ bits of $\Omega_V$
uniformly compute
$\DD_{n+f(n)-c}(W)$, for some constant $c$.}
\end{equation*}
In other words, the first $n$ bits of the halting probability of $V$ can
only solve the halting problem of $W$ for the inputs of length at most $n$
(plus or minus a constant).
We stress Tadaki's results apply to
arbitrary optimal machines -- not necessarily universal, and
characterize the computational relation between 
universal halting problems and universal halting probabilities
{\em for \pf machines}.

\subsection{Omega and \lce reals or c.e.\ sets}

Since Omega is Turing-complete, it computes all \lce reals and all \ce sets. We can then
obtain asymptotics regarding these computations. 
How many bits of Omega are needed in order
to compute $n$ bits of an arbitrary \lce real or an arbitrary \ce set?

In \cite{MR3521997} it was shown that given any computable $h:\Nat\to\Nat$,
\begin{equation}\label{tSjUiw9u4d}
\parbox{14cm}{if $\sum_n 2^{n-h(n)}<\infty$ converges, then Omega
computes every \lce real $\alpha$ with oracle-use $h$.}
\end{equation}
Moreover the converse  of this implication was shown under the additional assumption that
 $h(n)-n$ is non-decreasing: if the sum in \eqref{tSjUiw9u4d} is unbounded,
 then there are \lce reals which are not computable from Omega with oracle-use $h$.
 In fact, the following stronger statement was obtained:
\begin{equation*}
\parbox{14cm}{If $\sum_n 2^{n-h(n)}=\infty$  
then there exist two c.e.\  reals such 
that no \lce  real can compute both of them with oracle-use $h+\bigo{1}$.}
\end{equation*}
Hence we see that Omega computes all \lce reals with oracle-use $n+2\log n$ but 
oracle-use $n+\log n$ is not always sufficient for this purpose.
There is also a uniform (hence stronger) version of the latter fact, which was obtained recently by
Fang Nan:\footnote{A special case of this result for 
 oracle-use $n\mapsto n+\bigo{1}$ was obtained earlier in \cite{jflBaiasL06}.}
\begin{equation*}
\parbox{14cm}{There exists a \lce real which is not computed
by {\em any} Omega number with oracle-use 
$n\mapsto n+g(n)$, for any computable non-decreasing function $g$ such that 
$\sum_n 2^{-g(n)}=\infty$.}
\end{equation*}
Analogous results were obtained in \cite{MR3521997,Kobayashi_rep}
with respect to \ce sets, although now the oracle-uses are appropriately tighter.
Given any computable function $g$,
\begin{equation*}
\parbox{14cm}{if $\sum_n 2^{-g(n)}<\infty$, 
then every \ce set is computable from 
Omega with oracle-use $g$.}
\end{equation*}
Here, just as in \eqref{tSjUiw9u4d}, the statements hold for all versions of Omega.
A strong converse is also given, under the additional assumption that 
$g$ is nondecreasing:
\begin{equation*}
\parbox{14cm}{if $\sum_n 2^{-g(n)}=\infty$ 
then Omega  cannot compute all \ce sets
with oracle-use $g$. In fact, in this case, 
no linearly complete \ce set can be computed by any \ce real
with oracle-use $g$.}
\end{equation*}

One can draw interesting conclusions about the computational relation between
halting probabilities and halting problems from the last two results. Note that
the halting set $H$ with respect to plain Turing machines and a Kolmogorov numbering
of all programs, is linearly complete. Hence any such canonical halting problem $H$
is not computable by Omega with oracle-use $\log n$ but it is computable
by Omega with oracle-use $2\log n$. We may also contrast these results with
Tadaki's results that we discussed in Section \ref{8VvCMdWp4N}, which
referred to halting problems with respect to \pf machines.
 
These results give worse-case bounds on the number of bits of Omega needed for
the computations. One 
would guess that the number of bits of Omega that are needed for the computation of
$n$ bits of a \lce real or a \ce set $X$ should depend on how much information is encoded in
$X\restr_n$. The latter, in turn, would be reflected in the Kolmogorov complexity of the
first $n$ bits of $X$.
An upper bound on the oracle-use along these lines was obtained in 
\cite{koba_rod}:
\begin{equation}\label{tEeG4sH537}
\parbox{14cm}{Every \lce real $X$ can be computed from $\Omega$ with oracle-use
$g(n)=\min_{i\geq n} K(X\restr_n)$.}
\end{equation}
This is remarkable! The \pf complexity $K(X\restr_n)$ is supposed to measure the amount
of information -- in bits -- encoded in the first $n$ bits of $X$. Then
\eqref{tEeG4sH537} says that the number of bits of $\Omega$ 
containing the information $X\restr_n$  is $K(X\restr_n)$, \ie precisely 
the amount of information in $X\restr_n$; and $\Omega$ has this property with respect to
every \lce real $X$. This is another testimony of the compactness of retrievable information
in the initial segments of Omega.

\subsection{How similar or different are two Omega numbers?}\label{YfgHxGw4Ku}
In Kolmogorov complexity we know that changing universal machines is like
changing coordinate system in geometry: the theory remains the same, and
the complexity function only changes by a constant.
How is the change of the universal machine reflected in the halting probabilty?
How similar or how different are the halting probabilities with respect to two different
universal machines?

We know that all (universal) halting probabilities are in the same Turing degree, the degree
of the halting problem. So they can compute each other.
Moreover the oracle-use 
functions in these reductions are computably bounded, so all Omega numbers
are in the same weak truth-table degree as the halting problem \cite{Calude.Nies:nd}.
However these facts no longer hold for truth-table reductions. For example, the halting problem
is not truth-table reducible to any Omega number \cite{Calude.Nies:nd}.\footnote{This is
a consequence of a more general theorem of Demuth \cite{Demuth:88} 
which says that any non-computable 
set which is truth-table reducible to a \ml random set is Turing equivalent to a \ml random set.}
  
In \cite{jc/FigueiraSW06} it was shown that there are  universal $U,V$ such that
the corresponding halting probabilities have incomparable truth-table degrees.
In \cite{lata} it was shown that for each universal \pf machine $U$
there exists another universal \pf machine $V$ such that $\Omega_U\neq \Omega_V+\beta$
for every \lce real. This shows that given two 
Omega numbers, one is not necessarily the translation
of the other by a \lce real.

But how many bits
of $\Omega_U$ are needed to compute $\Omega_V$ given any 
two universal \pf machines $U,V$?
In \cite{MR3583619} it was shown that for each $\epsilon>1$ the oracle-use function 
$n\mapsto n+\epsilon\cdot\log n$ suffices for this computation while, in general, 
$n\mapsto n+\log n$ does not.\footnote{The actual results are more general. Given
any Omega number $\Omega$ there exists another Omega number which is not
computable from $\Omega$ with oracle-use $n+\log n$. Moreover, given any non-decreasing
computable function $g$ such that $\sum_i 2^{-g(i)}$ is finite, any Omega number is
computable from any other Omega number with oracle-use $n+n\mapsto g(n)$. An exact
characterization of the required worse-case redundancy is still lacking. For example,
is oracle-use $n+\log n+ 2\log\log n$ sufficient for the 
computation of any omega number from any other
Omega number?}

\subsection{Computational power versus randomness of Omega}\label{eGq9lCDoAE}
We have discussed the computational power of Omega in
Sections \ref{lNZ1nm3t28} and \ref{yKBcWTkPDT}. From these discussions it is clear that
Omega is a rather special algorithmically random sequence: it is Turing-complete\footnote{In fact,
the standard notion of \ml randomness that we use allows the existence of random numbers
that compute any given sequence. In particular, given any sequence $A$, by the \KG theorem,
there exists a \ml random $X$ which computes $A$.} and a \lce real.
On the other hand, intuitively one would expect that algorithmically random numbers
do not solve interesting problems such as the halting problem. Denis Hirschfeldt in many of his
expository talks distinguishes two kinds of algorithmically random numbers: the ones that are
`truly' random and the ones, such as Omega, 
that pass the \ml tests because they have enough information
so that they can `fake' randomness. In fact, this distinction is quantifiable and reflected in many
theorems. If one slightly increases the strength of randomness under consideration, then
all random numbers are necessarily incomplete. In other words, Omega has just enough randomness
to qualify as \ml random, but no more. {\em Difference randomness} \cite{FrNgDiff} is exactly
this slight strengthening of \ml randomness which characterizes the \ml random reals which do not compute the halting problem. Furthermore, incomplete \ml random reals not only fail to compute the halting problem, but are computationally weak in many other ways\footnote{For example, 
they do not compute
any complete extension of Peano arithmetic \cite{MR2258713frank}.} 
in line with our intuition about algorithmically random numbers.  

All these facts show that Omega is a computationally powerful random real, which is necessarily 
not {\em very} random (in the sense that it fails slightly stronger notions of algorithmic randomness).
In all of the examples from Sections \ref{lNZ1nm3t28} and \ref{yKBcWTkPDT} demonstrating the computational power of Omega, computations use entire initial segments of Omega rather than
individual bits. Is this necessary? In other words, can we still solve interesting problems if
we only have access to certain, possibly non-adjacent, bits of Omega?
This question has a strongly negative answer, for the following reasons.
If useful information was coded into individual bits of Omega inside
an infinite computable set $A$ with infinite complement, 
one would be able to solve interesting problems by truth-table queries on 
the sequence of bits $\Omega_A$ of Omega restricted on positions in $A$. 
However it is known that such a sequence $\Omega_A$ is Turing-incomplete and 
\ml random, so
by \cite{Demuth:88} any non-computable 
set which is truth-table reducible to $\Omega_A$ is Turing equivalent to an
incomplete \ml random set.
But we know that such sets
do not code any problems of interest, which are mostly located inside \ce degrees, and in fact
the complete \ce degree.
 
Using more advanced results, we can obtain further demonstrations that
parts of Omega such as its even bits do not compute useful non-trivial problems, or
even any considerably non-trivial problem.
Let $\Omega_0, \Omega_1$ be the even and odd bit-sequences of Omega respectively.
Then $\Omega_0$ is an incomplete 1-random real below $\mathbf{0}'$. Such reals are known
to have small computational power -- for example, any \ce set that they compute is 
$K$-trivial.
As a result, not much useful information can be recovered by querying  the even bits of Omega.
A similar statement is true for any nontrivial computable subset of bits of Omega, which means that
information is not coded into individual bits of Omega but rather into its initial segments.
Having said this, there are non-computable \ce sets that are computable from both $\Omega_0$
and $\Omega_1$. Moreover, a \ce set is computable from $\Omega_0$ if and only
if it is computable from $\Omega_1$ \cite{kpartsommore}. 
One can generalize this by partitioning the bits of Omega into $k$ many
disjoint families, and considering the sets that are computable by all of these $k$ many reals.
Interestingly, one then gets a strictly decreasing sequence of 
nontrivial subclasses of the $K$-trivial sets \cite{kpartsome}.

\section{Halting probability relative to a set}
There are more than one ways that one can relativize the halting probability --
here we consider two. In Section \ref{ITRpcDH23Q}
we discuss the probability $\Omega_U(X)$ that the universal \pf machine $U$ 
halts with an output
inside a given set $X$, focusing on the randomness properties of it.
In Section \ref{SSwIvk3kXz} we instead add the given set $X$ as an oracle to
the \pf machine $U$, and consider the halting probability $\Omega_U^X$ of
$U$ with oracle $X$.
Finally in Section \ref{eGq9lCDoAE} we report on some recent work concerning computations
from subsequences of the binary expansion of Omega, and connections with lowness classes
and $K$-triviality.

\subsection{Restricting the output of a universal \pf machine}\label{ITRpcDH23Q}
Given a universal \pf machine $U$ we may consider the probability 
$\Omega_U(X)$ that $U$ halts
with output in a given set $X$ of strings. Such probabilities were considered in
\cite{ipl/Kobayashi93,jsyml/BecherFGM06,jc/FigueiraSW06,omegax}, initially
as an attempt to produce concrete numbers which are more random than Chaitin's Omega.
In \cite{jsyml/BecherFGM06} it was shown that if $X$ is $\Sigma^0_n$-complete
or $\Pi^0_n$-complete
for any $n> 1$ then $\Omega_U(X)$ is 1-random.\footnote{The same result for
an apparently stronger notion of completeness was obtained earlier in \cite{ipl/Kobayashi93}.} 
On the other hand, in the same paper it was shown that 
if  $n> 1$ and $X$ is any $\Sigma^0_n$ or $\Pi^0_n$ set then 
$\Omega_U(X)$ is not $n$-random. 

The case $n=1$ is particularly interesting. 
Chaitin \cite{chaitin2004algorithmic} had already noticed that if $X$ is any non-empty 
$\Sigma^0_1$ set, 
then $\Omega_U(X)$ is 1-random and left-c.e., so it is just another Omega number. 
The case when $X$ is \pz is the basis of 
\cite[Question 8.10]{MR2248590} 
which was also discussed in 
\cite{DBLP:journals/jsyml/BecherG05,jsyml/BecherFGM06,jc/FigueiraSW06},
and
remained open until recently.
In 
\cite{jc/FigueiraSW06} an optimal (in terms of Kolmogorov complexity) 
but not universal \pf machine $U$ was constructed, 
and a \pz set $X$ such that  $\Omega_U(X)$ is not \ml random. 
The question whether $\Omega_U(X)$ is always an omega number when 
$U$ is universal and $X$ is a non-empty \pz set, 
required a deeper understanding of
omega numbers and their differences.
This is hardly surprising since, when $X$ is a \pz set the real $\Omega_U(X)$
is the difference of two omega numbers.
In \cite{omegax} it was shown that $\Omega_U(X)$ is 
\ml random and \lce \ie just another omega number when $U$ is universal and $X$ is
any \pz set.
 The proof uses the results about differences of omega numbers that we
discussed in Section \ref{8kFnCcyiiP} 
and a decanter argument for the construction of a suitable
\ml test.

Many questions which ask for the randomness strength of $\Omega_U(X)$
when $X$ has some type of universality remain open. For example,
if $X$ is uniformly $\Sigma^0_n$-complete for all $n$ (\ie it is in the many-one degree of 
$\emptyset^{(\omega)}$) then what level of randomness can we expect from 
$\Omega_U(X)$?

\subsection{Halting probability in an oracle \pf machine}\label{SSwIvk3kXz}
Another way to relativize the halting probability is to consider
universal oracle \pf machines $U$ and the probability 
$\Omega_U^X$ that they halt when equiped with oracle $X$.
The study of these numbers, as well as the map $X\mapsto \Omega_U^X$
was initiated in
\cite{Downey.Hirschfeldt.ea:05}.
First, we note that the Omega characterization \eqref{nIfPvZ8kto} 
relativizes to any oracle $X$ in 
a rather straightforward way. For example, the halting probabilities of  universal
\pf machines with oracle $X$ are exactly the $X$-left-c.e.\  real numbers which are
1-random relative to $X$.
In the same way, Solovay reducibility relativizes to $X$ and the numbers $\Omega^X$
are exactly the $X$-left-c.e.\  real numbers with the slowest $X$-left-c.e.\  approximations.
Second, in \cite{Downey.Hirschfeldt.ea:05} it was shown that
\[
\parbox{14cm}{$\Omega$ is 1-random relative to $X$ iff $\Omega^X_U$ is \lce for some
universal \pf machine $U$.}
\]
Sets $X$ with this property are called {\em low for $\Omega$}.
Recall that $X$ is $K$-trivial if
there exists some $c$ such that $K(X\restr_n)\leq K(n)+c$ for all $n$,
and note that in general $X\not\leq_T \Omega^X_U$.
a remarkable result from \cite{Downey.Hirschfeldt.ea:05} is that
\[
\parbox{14cm}{$X$ is $K$-trivial iff $\Omega^X_U$ is \lce for all
universal \pf machines $U$.}
\]
and moreover, $X$ is $K$-trivial if and only if 
$X'\equiv_T \Omega^X_U$.

Perhaps the leading motivating question behind \cite{Downey.Hirschfeldt.ea:05}
was the degree invariance of the operator $X\mapsto \Omega_U^X$. A very strong
negative answer was given to this question: 
\[
\parbox{14cm}{given any universal
\pf machine $U$,  there are sets $A,B$ which only differ on finitely many bits,
and such that $\Omega_U^A$ is Turing incomparable to
$\Omega_U^B$.}
\]
Another very important consequence of this work concerns the invariance of
the relativized Omega with respect to the choice of the underlying universal \pf machine.
Recall that Chaitin's omega with respect to different universal \pf machines are
very similar -- certainly Turing equivalent. Is this true for the relativized Omega?
If the oracle $X$ is $K$-trivial, then it follows from the above discussion that
$\Omega_U^X$ is always in the Turing degree of the halting problem. A remarkable fact
is that 
\[
\parbox{14cm}{if $X$ is not $K$-trivial, then there are universal Turing machines $U,V$ such that
$\Omega_U^X$ and $\Omega_V^X$ have incomparable Turing degree.}
\]
In \cite{joelowklowwea,MR2835602} it was shown that the low for $\Omega$ reals
are very related to relativized \pf complexity. In particular,
the following are equivalent for each $X$:
\begin{enumerate}[\hspace{0.5cm}(a)]
\item $X$ is low for $\Omega$;
\item there exists $c$ such that $K(n)\leq K^X(n)+c$ for infinitely many $n$;
\item $\Big\{Z\ |\ \forall n\ K^X(n)\leq K^Z(n)+c\Big\}$ is finite for each $c$;
\item $\Big\{Z\ |\ \forall n\ K^X(n)\leq K^Z(n)+c\Big\}$ is countable for each $c$.
\end{enumerate}

Finally low for $\Omega$ sets $X$ 
have interesting growth-rate properties, as we briefly mentioned
in Section \ref{QgvOh5fZzw}. 
If $X$ is non-computable and low for $\Omega$ then it computes a function which is not
dominated by any computable function (by Miller and Nies, see \cite[Section 8.1]{Ottobook}).
On the other hand, by \cite{MR3248779} if $X$ is low for $\Omega$ then any
function $f$ which is computable by $X$ does not  speed-up the effective approximations to
Omega, in the sense that there exists a constant $c$ such that 
$\Omega-\Omega_n\leq c\cdot (\Omega-\Omega_{f(n)})$ for all $n$.

\subsection{Omega operators}
Omega can be turned into an operator, mapping binary streams or subsets of $\Nat$ to reals, in
a number of different ways. Perhaps the more natural such examples are the following:
%
\begin{eqnarray*}
\textrm{(a)}\hspace{0.3cm} Z\mapsto \sum_{U(\sigma)\de\in Z} 2^{-|\sigma|}\hspace{1cm} 
&&\textrm{(b)}\hspace{0.3cm}  X\mapsto \sum_{U^X(\sigma)\de} 2^{-|\sigma|}\\[0.5cm]
\textrm{(c)}\hspace{0.3cm} X\mapsto \sum_{n} 2^{-K(X\restr_n)}\hspace{1cm}
&&\textrm{(d)}\hspace{0.3cm}  X\mapsto \sum_{U(\sigma)\de\prec X} 2^{-|\sigma|}
\end{eqnarray*}
where $Z\subseteq\Nat$, $X$ is a binary stream and `$\prec$' denotes the prefix relation.

The main motivation for the study of operator (a) in \cite{DBLP:journals/jsyml/BecherG05,jsyml/BecherFGM06, tcs/BecherG07,DBLP:journals/jsyml/BecherG09} was
the discovery of concrete highly random numbers. Some facts about operator (a)
were discussed in Section \ref{ITRpcDH23Q}.
Much about the initial excitement regarding operator (b)
had to do with the possibility that it provides a uniform degree-invariant 
solution to Post's problem, which is
a long-standing question in classical computability theory. 
As we discussed in Section \ref{SSwIvk3kXz},
the results in \cite{Downey.Hirschfeldt.ea:05} provide a strong negative answer to such a hope.

Several results regarding the analytic behavior of operator (b), focused on its continuity and
the complexity of its range, were obtained 
in  \cite[Section 9]{Downey.Hirschfeldt.ea:05}. It was shown that 
%
\begin{itemize}
\item Operator (b) is lower semi-continuous but not continuous;
\item Operator (b) is continuous exactly on the 1-generic reals.
\end{itemize}
Remarkably, the supremum of (b) is achieved, and the argument $X$ that achieves it is always
1-generic. Several questions about the nature of operator (b) remain, and the 
interested reader is referred to \cite{Downey.Hirschfeldt.ea:05}.

Recently a study of operator (c) was initiated in \cite{HMSY}, where it was shown that

\begin{itemize}
\item operator (c) is continuous and almost everywhere differentiable with derivative 0;
\item operator (c) is co-meagerly non-differentiable and nowhere monotonic.
\end{itemize}
In fact, operator (c) is differentiable exactly on the \ml random reals. 
Note that, strictly speaking, the value of (c) on any real depends on the choice of the underlying
universal machine $U$. 
One can construct a universal machine $U$ with respect to which 
operator (c) maps any random stream to a non-random real.
On the other hand, 
there exists a stream $X$ which is always mapped to a non-random real via (c), independently of
the underlying optimal machine $U$.
It was also shown in \cite{HMSY}  that the inverse image of (c) on any real is null, for a certain
underlying optimal machine; it is an open question whether this holds for all optimal underlying machines.
Finally
we note that operator (d) is very similar to operator (c), and  the results we presented about (c) 
(all from \cite{HMSY}) also hold for (d) by the same or similar proofs.

\section{Machine probabilities beyond halting}
In turns out that the probability of many properties of a universal randomized
Turing machine can be expressed in terms of relativized halting probabilities.
This phenomenon was the topic of \cite{MR3629268}.
In Section \ref{UnEmmQPfiV} we give such examples and their characterizations in 
terms of algorithmic randomness. 
Then in Section \ref{PN5XnbpwfZ} we note that although there is usually a rough correspondence
between the arithmetical complexity of a property and the strength of the algorithmic randomness
of its probability, this is no always the case, even for some very natural properties such as
a computable infinite output from a randomized universal machine. 
Finally in Section \ref{7dkEgGlGGy} we note that for properties of higher arithmetical complexity
than halting (for example 2-quantifier definable properties) the degree or even the
algorithmic randomness  or computability of the probability is very sensitive to
the choice of underlying universal
machine, in contrast with the halting probability. 

\subsection{Natural machine properties and their universal probabilities}\label{UnEmmQPfiV}
For example, consider an oracle Turing machine $M$ which has one tape for the oracle $X$, 
one input tape for the input string $\sigma$ and one output tape for the output
string, in case the computation halts. Then we may treat the oracle $X$ 
as a random variable, so that the machine is viewed as a probabilistic machine. Then we
can consider the probability that the partial function $\sigma\mapsto M^X(\sigma)$ 
has certain properties, when $M$ is universal. 
The probability that $\sigma\mapsto M^X(\sigma)$ is total (the {\em totality probability} of $M$)
is the same as the non-halting probability of a universal \pf machine with oracle the halting problem
$\emptyset'$. In other words, given a universal oracle Turing machine $M$, there exists a
universal oracle \pf machine $U$ such that  the totality probability of $M$ is equal to 
$1-\Omega_U^{\emptyset'}$; vice-versa, given a universal oracle \pf machine $U$
there exists a universal oracle Turing machine $M$ with totality probability equal to 
$1-\Omega_U^{\emptyset'}$. We depict this fact in Table \ref{GGvIp7p1CY} by
saying that totality (of a universal oracle machine) corresponds to $1-\Omega^{\emptyset'}$;
in other words, totality probability of a universal  oracle machine is a \rcep 2-random real.

Next, consider a universal oracle-machine $M$ which, given an oracle $X$, it enumerates
a set of strings $W^X$. Note that this can be formalized as a property of the
function $\sigma\mapsto M^X(\sigma)$, for example by letting  $W^X$ be the domain
(or the range) of this function. Again, if we consider $X$ as a random variable, we can
consider the probability that $W^X$ is a computable set.
It turns out that this probability
is the same as the 
halting probability of some
universal \pf machine with oracle $\emptyset''$.
We note this fact in Table \ref{GGvIp7p1CY} by saying that the property of
probabilistically
enumerating a computable set via a universal oracle Turing machine corresponds to
 $\Omega^{\emptyset^{(2)}}$.
We may continue with many of the properties that one meets in a first course in computability theory
such as the property that $W^X$ is co-finite or the property that 
$W^X$ computes the halting problem. These characterizations are included in
Table \ref{GGvIp7p1CY}.
Similar characterizations may be obtained for different models of Turing machines
even when there are no oracles present. For example, given a \pf machine $U$,
Wallace (see \cite[Section 0.2.2]{Dowe2008a} and \cite[Section 2.5]{Dowe2011a})
considered the measure of reals $X$ such that the \pf machine
$\sigma\mapsto U(X\restr_n\ast\sigma)$ is universal for all $n$, and called it the
{\em universality probability of $U$}. In \cite{Barmpalias3488} it was shown
that the universality probabilities of universal \pf machines are exactly the
non-halting probabilities relative to $\emptyset^{(3)}$.

\subsection{Complexity of a property versus
algorithmic randomness of its probability}\label{PN5XnbpwfZ}
The reader may have noticed a pattern in the examples of Table \ref{GGvIp7p1CY}. 
The arithmetical complexity
of each of the properties we considered matches the level of randomness
that they correspond to.  For example, totality is a $\Pi^0_2$ property and
it corresponds to 2-random \rcep reals, which is the analogue of $\Pi^0_2$ in
algorithmically random reals.  Indeed, recall that
the \rcep reals are exactly the reals whose left Dedekind cut is $\Pi^0_2$.
Similarly, computability is a $\Sigma^0_3$ property\footnote{For example, the set
$\{e\ |\ W_e \ \textrm{is computable}\}$ is $\Sigma^0_3$.}
and the property that $W^X$ is computable corresponds to the 3-random \lcepp reals
$\Omega^{\emptyset^{(2)}}$ which are the analogue of $\Sigma^0_3$ in
algorithmically random reals. Finally, universality is a $\Pi^0_4$ property and, sure enough,
the universality probabilities are characterized by the
numbers $1-\Omega^{\emptyset^{(3)}}$ which are the analogues of $\Pi^0_4$ for
algorithmically random reals.

\begin{table}
\centering
\begin{tabular}{llc}
 \hline\cmidrule{1-3} 
 {\small Totality} & &   {\small $1-\Omega^{\emptyset'}$}\\[0.1cm]
 {\small Enumeration of a computable set}	&   &   {\small$\Omega^{\emptyset^{(2)}}$}\\[0.1cm]
{\small  Enumeration of a co-finite set}   &  & {\small$\Omega^{\emptyset^{(2)}}$}\\[0.1cm]
{\small  Enumeration of a set which computes 
$\emptyset'$} &  &  {\small$\Omega^{\emptyset^{(3)}}$}\\[0.1cm]
{\small  Universality probability}   & & {\small$1-\Omega^{\emptyset^{(3)}}$}\\[0.1cm]
 \hline\cmidrule{1-3} 
\end{tabular}
\caption{Universal probabilities as Omega numbers}
\label{GGvIp7p1CY}
\end{table}

Is this a coincidence, or is there a general theorem that characterizes the various universal
probabilities in terms of algorithmically random numbers of the same arithmetical complexity?
With so many examples and a uniform methodology for their analysis (see \cite{MR3629268})
it is tempting to guess a positive answer to this question.
Quite surprisingly, there is a simple example of a universal probability which
behaves very differently and, in fact, its complexity depends strongly on the
choice of the universal machine. Given a universal oracle Turing machine $M$ consider the
probability that $n\mapsto M^X(n)$ is computable (and hence, total). This question can also be
phrased in terms of monotone machines, where the input is a stream $X$ and the output is
either a string or a stream $Y$, and we are looking for the probability that the output
is a computable stream. Note that the property in question is $\Sigma^0_3$, so according to
the previous discussion, one would expect these probabilities to be characterized as 
$\Omega^{\emptyset^{(2)}}$. However it was shown in \cite{procom} that this is never a 
3-random number, and depending on the underlying universal machine it can be
as complex as $\Omega^{\emptyset'}$ or $1-\Omega^{\emptyset'}$,
and as simple as 1/2.

\subsection{Invariance with respect to different universal machines}\label{7dkEgGlGGy}
We have seen that, although certain differences may exist between two universal halting probabilities,
their Turing degree remains fixed to the degree of the halting problem. This is no longer true
for probabilities of properties of higher complexity, such as the ones we
discussed in Section \ref{UnEmmQPfiV}. In many cases,
if one examines the proofs of this failure of invariance, one reason is the fact from 
\cite{Downey.Hirschfeldt.ea:05} that the Turing degree of relativizations of Omega are very sensitive to the choice of the universal \pf machine.\footnote{Of course the probabilities 
of Section \ref{UnEmmQPfiV}  are not with respect to oracle machines,
but as Table \ref{GGvIp7p1CY} shows they can be characterized as such.}
For example, in \cite{Barmpalias3488} it was shown that there are two universal \pf machines
with universality probabilities having incomparable Turing degrees. Similar results were obtained in
\cite{MR3629268} with regard to the other probabilities that we discussed. 

Perhaps the most extreme example with respect to sensitivity 
over the choice of underlying universal machine $U$ is
the probability of producing an infinite computable stream -- the example discussed in
Section \ref{PN5XnbpwfZ}. In this case it is not only the Turing degree of the probability that
depends on $U$, but most computational aspects of this number. Depending on the choice of $U$
it could be as simple as 1/2 or as complicated as $\Omega^{\emptyset'}$ or even the mirror image
 $1-\Omega^{\emptyset'}$ (namely a \rce real relative to the halting problem $\emptyset'$ which is
 \ml random relative to $\emptyset'$).

\section{Conclusion}
We have discussed several mathematical aspects of Chaitin's halting probability, 
which go well beyond its simple algorithmic properties that were established Chaitin's seminal paper.
We have also provided a rich bibliography on this topic, where further results and properties may 
be found, which have not found a place in our discussions. We hope we have shown that
Chaitin's Omega plays an important role in contemporary studies in algorithmic randomness.


\end{document}